\documentclass{amsart}

\newif\ifdraft\draftfalse
\newif\ifcite\citefalse
\newif\ifblow\blowtrue

\usepackage{amsfonts,amssymb, amsmath}
\usepackage{oldgerm}
\usepackage{amscd}
\usepackage{eucal} 
\usepackage{mathrsfs} 
\usepackage[hyperindex]{hyperref}  
\usepackage[alphabetic]{amsrefs}
\ifdraft\ifcite\usepackage{showkeys}\else\usepackage[notcite,notref]{showkeys}\fi\fi

\usepackage{xypic}
\newdir{ >}{{}*!/-5pt/\dir{>}}

\newtheorem{proposition}[equation]{Proposition}
\newtheorem{theorem}[equation]{Theorem}

\theoremstyle{remark}

\theoremstyle{definition}
\newtheorem{definition}[equation]{Definition}

\newtheorem{para}[equation]{}

\theoremstyle{remark}

\newtheorem{remark}[equation]{Remark}

\numberwithin{equation}{section}

\def\bz{\mathbb Z}

\def\bp{\mathbb P}

\def\bean{\begin{eqnarray}}
\def\eean{\end{eqnarray}}
\def\bea{\begin{eqnarray*}}
\def\eea{\end{eqnarray*}}

\def\x0{X_{s_0}}
\def\2p{\bp^1\times \bp^1}

\def\a{\alpha}

\def\nm{\nonumber}
\def\er{\eqref}

\def\vol{d\sigma}

\def\ben{\begin{equation}}
\def\een{\end{equation}}

\newcommand\calL{\mathcal L}

\newcommand\wt{\widetilde}

\newcommand\ind{\mathop{\mathrm{Ind}}\nolimits}
\newcommand\ppv{\partial_+ V}
\newcommand\pmv{\partial_- V}
\newcommand\vn{\vec n}
\newcommand\sing{\mathop{\mathrm{Sing}}\nolimits}

\newcommand\om{\omega}
\newcommand\Om{\Omega}
\newcommand\Ome{\Omega^M}
\newcommand\so{\mathrm{SO}}
\newcommand\ep{\epsilon}
\newcommand\dph{\frac{\partial}{\partial \phi}}
\newcommand\up{\Upsilon}

\begin{document}
\title{Secondary Chern-Euler forms and the Law of Vector Fields}
\author{Zhaohu Nie}
\email{znie@psu.edu}
\address{Department of Mathematics\\
Penn State Altoona\\
3000 Ivyside Park\\
Altoona, PA 16601, USA}

\date{\today}

\begin{abstract}
The Law of Vector Fields is a term coined by Gottlieb for a relative Poincar\'e-Hopf theorem. It was first proved by Morse \cite{morse} and expresses the Euler characteristic of a manifold with boundary in terms of the indices of a generic vector field and the inner part of its tangential projection on the boundary. We give two elementary differential-geometric proofs of this topological theorem, in which secondary Chern-Euler forms  \cite{chern} naturally play an essential role. In the first proof, the main point is to construct a chain away from some singularities. The second proof employs a  study of 
the secondary Chern-Euler form on the boundary, which may be of independent interest. More precisely, we show by explicitly constructing a primitive that, away from the outward and inward unit normal vectors, the secondary Chern-Euler form is exact up to a pullback form. In either case, Stokes' theorem is used to complete the proof.
\end{abstract}

\maketitle

\section{Introduction}

Let $X$ be a smooth oriented compact Riemannian manifold with boundary $M$. Throughout the paper we fix $\dim X=n\geq 2$ and hence $\dim M=n-1$. 
On $M$, we have a canonical decomposition
\begin{equation}
\label{collar}
TX|_M\cong \nu\oplus TM,
\end{equation}
where $\nu$ is the rank 1 trivial normal bundle of $M$. 

Let $V$ be a smooth vector field on $X$.
We assume that $V$ has only isolated singularities, i.e., the set $\sing V:=\{x\in X|V(x)=0\}$ is finite, and that the restriction $V|_M$ is nowhere zero. Define the index $\ind_x V$ of $V$ at an isolated singularity $x$ as usual (see, e.g., \cite[p. 136]{hirsch}), and let $\ind V=\sum_{x\in \sing V} \ind_x V$ 
denote the sum of the local indices. 

\begin{para}\label{n&N}
As an important special case, let $\vn$ be the outward unit normal vector field of $M$, and $\vec N$ a generic extension of $\vn$ to $X$. Then by definition 
\ben
\label{standard0} 
\ind \vec N=\chi(X),
\een
where $\chi(X)$ is the Euler characteristic of $X$ (see, e.g., \cite[p. 135]{hirsch}). 
\end{para}

For a general $V$, let $\partial V$ be the projection of $V|_M$ to $TM$ according to \er{collar}, and let 
$\partial_- V$ (resp. $\partial_+ V$) be the restriction of $\partial V$ to the subspace of $M$ where $V$ points inward (resp. outward) to $X$. Generically $\partial_\pm V$ have isolated singularities. (A non-generic $V$ can always be modified by adding an extension to $X$ of a normal vector field or a tangent vector field to $M$.)

Using 
the flow along $-V$ and counting fixed points with multiplicities, we have the following 
\emph{Law of Vector Fields}:
\ben
\label{law of v f}
\ind V+\ind \pmv=\chi(X).
\een
Naturally this is a relative Poincar\'e-Hopf theorem. It was first proved by Morse \cite{morse} and later on publicized by Gottlieb, who also coined the name ``Law of Vector Fields". 

One main purpose of this paper is to give two elementary differential-geometric proofs of this theorem \er{law of v f}. 

In his famous proof \cite{chern} of the Gauss-Bonnet theorem, Chern constructed a differential form $\Phi$ (see \er{phi}) of degree $n-1$ 
on the tangent sphere bundle $STX$, consisting of unit vectors in $TX$, satisfying the following two conditions:
\begin{equation*}
\label{dP=O}
d\Phi=-\Omega,
\end{equation*}
where $\Omega$ is the Euler curvature form of $X$ (pulled back to $STX$), which is defined to be $0$ when $\dim X$ is odd (see \er{omega}), and 
\begin{equation*}
\wt{\Phi_0}=\wt{\vol}_{n-1},
\end{equation*}
i.e., the 0th term $\wt{\Phi_0}$ of $\Phi$ is the relative unit volume form for the fibration $S^{n-1}\to STX \to X$
(see \er{0th term}). We call $\Phi$ the \emph{secondary Chern-Euler form}. 

Define $\alpha_V:M\to STX|_M$ by rescaling $V$, i.e., $\a_V(x)=\frac{V(x)}{|V(x)|}$ for $x\in M$. Then Chern's basic method \cite[\S 2]{chern2}\cite[\S 6]{bc} 
using the above two conditions
and Stokes' theorem gives 
\ben
\label{general}
\int_X \Omega=-\int_{\alpha_V(M)}\Phi+\ind V
\een
(see \er{basic method}). 
Applying \er{general} to the $\vec n$ and $\vec N$ in \ref{n&N} and using \er{standard0}, one gets the following relative Gauss-Bonnet theorem in \cite{chern2}
\ben
\label{standard}
\int_X \Omega=-\int_{\vn(M)}\Phi+\ind(\vec N)=-\int_{\vn(M)}\Phi+\chi(X).
\een
Comparison of \er{general} and \er{standard} gives
\ben
\label{from chern's method}
\chi(X)=\ind V+\int_{\vn(M)}\Phi-\int_{\alpha_V(M)}\Phi.
\een

The following is our main result that identifies \er{from chern's method} with the Law of Vector Fields \er{law of v f}. 
\begin{theorem}\label{two terms} 
The following formula holds:
$$
\int_{\vn(M)}\Phi-\int_{\alpha_V(M)}\Phi=\ind \pmv.
$$
\end{theorem}

A first proof of the above theorem is given in Section \ref{first proof}. The main point of this first proof is to construct, away from some singularities, a chain connecting $\a_V(M)$ to $\vec n(M)$ and then to apply the Stokes' theorem. 

A second proof of Theorem \ref{two terms} to be given in Section \ref{trans} employs a  study of the secondary Chern-Euler form on the boundary, i.e., when the structure group is reduced from $\so(n)$ to $1\times \so(n-1)$. This study may be of some independent interest. 

In more detail, the images $\vec n(M)$ and $(-\vec n)(M)$ in $STX|_M$ 
are the spaces of outward and inward unit normal vectors of $M$. Define
\begin{equation}
\label{eq-def-cstm}
CSTM:=STX|_M\backslash(\vec n(M)\cup (-\vec n)(M))
\end{equation}
 ($C$ for cylinder) to be the complement. Also let $\pi:STX|_M\to M$ be the natural projection. 

\begin{theorem}\label{new main} There exists a differential form $\Gamma$ of degree $n-2$ on CSTM, such that after restricting to $CSTM$ 
\ben
\label{main formula}
\Phi-\pi^*\vec n^* \Phi=d\Gamma.
\een
\end{theorem}

The form $\Gamma$ is defined in \er{def Gamma}, and the above theorem is proved right after that utilizing Propositions \ref{much quicker} and \ref{big Phi}. 

At the end of Section \ref{trans}, we employ Stokes' theorem to give a second proof of Theorem \ref{two terms}, and hence of the Law of 
Vector Fields \er{law of v f}, using Theorem \ref{new main}.

\begin{remark} Unlike in \cite{sha} or \cite{nie}, we do not assume that the metric on $X$ is locally product near its boundary $M$. Therefore our results in this paper deal with the general case and generalize those in \cite{nie}. 
\end{remark}

\begin{remark} We would like to emphasize the elementary nature of our approaches, in the classical spirit of Chern in \cites{chern,chern2}. Transgression of Euler classes has gone through some modern development utilizing Berezin integrals. 
The Thom class in a vector bundle and its transgression are studied in \cite{MQ}. This Mathai-Quillen form is further studied in \cite{BZ} and \cite{BM}.
For the modern developments, we refer the reader to the above references and two books, \cite{BGV} and \cite{Z}, on this subject. 
\end{remark}

\section{Secondary Chern-Euler forms}

In this section, we review the construction, properties and usage of the secondary Chern-Euler form $\Phi$ in \cite{chern}, which plays an essential role in our approaches. 

Throughout the paper, $c_{r-1}$ denotes the volume of the unit $(r-1)$-sphere $S^{r-1}$. We also agree on the following ranges of indices
\ben
\label{range}
1\leq A,B\leq n,\ 2\leq \a,\beta\leq n-1,\ 2\leq s,t \leq n.
\een

The secondary Chern-Euler form $\Phi$ is defined as follows. 
Choose oriented local orthonormal frames $\{e_1,e_2,\cdots,e_n\}$ for the tangent bundle $TX$. Let $(\omega_{AB})$ and 
$(\Omega_{AB})$ be the $\mathfrak{so}(n)$-valued connection forms and curvature forms
for the Levi-Civita connection $\nabla$ of the Riemannian metric on $X$ defined by 
\begin{gather}
\nabla e_A=\sum_{B=1}^n \omega_{AB}e_B,\label{def-omega}\\
\Omega_{AB}=d\omega_{AB}-\sum_{C=1}^n \omega_{AC}\omega_{CB}.\label{def-O}
\end{gather}
(In this paper, products of  differential forms always mean ``exterior products" although we omit the notation $\wedge$ for simplicity.) 
Let the $u_A$ be the coordinate functions  on $STX$ in terms of the frames defined by
\ben
\label{def-u}
v=\sum_{A=1}^{n} u_A(v)e_A,\quad \forall v\in STX.
\een
Let the $\theta_A$ be the 1-forms on $STX$ defined by
\ben
\label{def-theta}
\theta_A=du_A+\sum_{B=1}^{n}u_B\omega_{BA}.
\een
For $k=0,1,\cdots,[\frac{n-1}2]$ (with $[-]$ standing for the integral part), define degree $n-1$ forms on $STX$
\ben
\label{eq-phi-j}
\Phi_k=\sum_{A} \epsilon(A)u_{A_1}\theta_{A_2}\cdots\theta_{A_{n-2k}}\Omega_{A_{n-2k+1}A_{n-2k+2}}\cdots \Omega_{A_{n-1}A_{n}},
\een
where the summation runs over all permutations $A$ of $\{1,2,\cdots,n\}$, and $\epsilon(A)$ is the sign of $A$. (The index $k$ stands for the number of curvature forms involved. Hence the restriction $0\leq k\leq [\frac{n-1}2]$. This convention applies throughout the paper.) Define the secondary Chern-Euler form as
\begin{align}
\Phi=\frac{1}{(n-2)!!c_{n-1}} \sum_{k=0}^{[\frac{n-1}2]} (-1)^k \frac{1}{2^k k! (n-2k-1)!!} \Phi_k\label{phi}
=:\sum_{k=0}^{[\frac{n-1}2]} \wt{\Phi_k}.
\end{align}
The $\Phi_k$ and hence $\Phi$ are invariant under $\mathrm{SO}(n)$-transformations of the local frames and hence are intrinsically defined. Note that the 0th term  
\begin{equation}
\label{0th term}
\wt{\Phi_0}=\frac{1}{(n-2)!!c_{n-1}} \frac 1{(n-1)!!}\Phi_0=\frac{1}{c_{n-1}}d\sigma_{n-1}=\wt{\vol}_{n-1}
\end{equation}
is the relative unit volume form of the fibration $S^{n-1}\to STX\to X$, since by \er{eq-phi-j}
\ben\label{factorial}
\Phi_0=\sum_A \epsilon(A)u_{A_1}\theta_{A_2}\cdots\theta_{A_n}=(n-1)!d\sigma_{n-1}
\een
(see \cite[(26)]{chern}). 

Then \cite[(23)]{chern} and \cite[(11)]{chern2} prove that 
\ben
\label{dphi}
d\Phi=-\Omega,
\een
where 
\ben
\label{omega}
\Omega=
\begin{cases}
0, & \text{if }n\text{ is odd,}\\
(-1)^{m}\frac{1}{(2\pi)^m 2^m m!}\sum_{A} \epsilon(A)\Omega_{A_1A_2}\cdots\Omega_{A_{n-1}A_n}, & \text{if }n=2m \text{ is even}
\end{cases}
\een
is the Euler curvature form of $X$. 

Now we review Chern's basic method \cite[\S 2]{chern2}\cite[\S 6]{bc}  of relating indices, $\Phi$ and $\Omega$ using Stokes' theorem. Similar procedures will be employed twice later. Let $V$ be a generic vector field on $X$ with isolated singularities $\sing V$. Let $B_r^X(\sing V)$ (resp. $S_r^X(\sing V))$ denote the union of small open balls (resp. spheres) of radii $r$ in $X$ around the finite set of points $\sing V$. Define $\alpha_V:X\backslash B_r^X(\sing V)\to STX$ by rescaling $V$. Then using \er{dphi} and Stokes' theorem, 
one proves \er{general} as 
\begin{align}
\label{basic method}
&\int_{X} \Omega=\lim_{r\to 0}\int_{\alpha_V(X-B_r^X(\sing V))} \Omega=\lim_{r\to 0}\int_{\alpha_V(X-B_r^X(\sing V))} -d\Phi\\
=&-\int_{\alpha_V(M)}\Phi+\lim_{r\to 0} \int_{\a_V(S^X_r(\sing V))} \Phi=-\int_{\alpha_V(M)}\Phi+\ind V\nm,
\end{align}
where the last equality follows from the definition of index and \er{0th term}. 

\section{First proof by constructing a chain}
\label{first proof}

In this section, we give a first proof of Theorem \ref{two terms} by constructing a chain, away from $\sing \pmv$, connecting $\alpha_V(M)$ to $\vec n(M)$. 

\begin{proof}[First proof of  Theorem \ref{two terms}] 
By definition, $\sing\pmv$ consists of a finite number of points $x\in M$ such that $\a_V(x)= -\vn(x)$. For $x\notin \sing\pmv$, let $C_x$ be the unique directed shortest great circle segment pointing from $\a_V(x)$ to $\vn(x)$ in $ST_xX$. With the obvious notation from before, let $U_r=M\backslash B_r^M(\sing\pmv)$ denote the complement in $M$ of the union of open balls of radii $r$ in $M$ around $\sing\pmv$. Obviously its boundary $\partial U_r=-S_r^M(\sing\pmv)$. 
Then
 \ben
 \label{boundary}
 \partial\left({\cup_{x\in U_r}C_x}\right)=\cup_{x\in U_r}\partial C_x-\cup_{x\in \partial U_r} C_x=
\vn(U_r)-\a_V(U_r)+W_r,
 \een
 with 
 \ben
 \label{wr}
 W_r:=\cup_{x\in S_r^M(\sing\pmv)}C_x. 
 \een
 Note the negative sign from graded differentiation in the second expression.
From \er{dphi} and \er{omega},  we have 
\ben
\label{d=0}
d\Phi=0\text{ on }STX|_M,
\een
since even if $\dim  X$ is even, $\Omega|_M=0$ by dimensional reason. \er{boundary}, Stokes' theorem and \er{d=0} imply
\begin{align}
\label{start}
&\int_{\vn(M)}\Phi-\int_{\alpha_V(M)}\Phi=\lim_{r\to 0} \left(\int_{\vn(U_r)}\Phi-\int_{\alpha_V(U_r)}\Phi\right)\\
=&-\lim_{r\to 0}\int_{W_r} \Phi=-\lim_{r\to 0}\int_{W_r} \wt{\Phi_0},\nm
\end{align}
where the last equality follows from \er{phi} and $\displaystyle\lim_{r\to 0}\int_{W_r} \wt{\Phi_k}=0$ for $k\geq 1$,  since such $\wt{\Phi_k}$'s in \er{eq-phi-j} involve curvature forms and don't contribute in the limit (see \cite[\S 2]{chern2}). 
By \er{0th term}, $\wt{\Phi_0}=\frac{1}{c_{n-1}}d\sigma_{n-1}$ is the relative unit volume form.
We then compute the RHS of \er{start} using spherical coordinates. 

\begin{para}\label{setup}
At $TX|_M$, we choose oriented local orthonormal frames $\{e_1,e_2,\cdots,e_n\}$ such that $e_1=\vec n$ is the outward unit normal vector of $M$. Therefore $(e_2,\cdots,e_{n})$ are oriented local orthonormal frames for $TM$. 
Let $\phi$ be the angle coordinate on $STX|_M$ defined by
\ben
\label{def-phi}
\phi(v)=\angle(v,e_1)=\angle(v,\vn),\ \forall v\in STX|_M. 
\een
we have from \er{def-u}
\begin{equation}
\label{eq-u1}
u_1=\cos\phi.
\end{equation}
Let
\begin{align}
p:CSTM=STX|_M\backslash(\vn(M)\cup (-\vn)(M))\to STM;\ v&\mapsto \frac{\partial v}{|\partial v|}
\label{def p}\\
\text{(in coordinates)} (\cos\phi,u_2,\cdots,u_{n})&\mapsto\frac{1}{\sin\phi}(u_2,\cdots,u_{n})\nm
\end{align}
be the projection to the equator $STM$.
By definition, 
\ben\label{pa}
p\circ \a_V=\a_{\partial V}\text{ when }\partial V\neq 0. 
\een
\end{para}

Therefore the image of $W_r$ in \er{wr} under the above projection is 
$$
p(W_r)=\cup_{x\in S_r^M(\sing\pmv)}p(C_x)=\cup_{x\in S_r^M(\sing\pmv)}\a_{\partial V}(x)=\a_{\partial V}\bigl(S_r^M(\sing\pmv)\bigr).
$$ 
On $C_x$ for $x\in M$, the $\phi$ \er{def-phi} ranges from $\phi(\alpha_V(x))$ to 0. 

The relative volume forms $d\sigma_{n-1}$ of $S^{n-1}\to STX|_M\to M$ and $d\sigma_{n-2}$ of $S^{n-2}\to STM\to M$ are related by 
\ben
\label{vol relation}
d\sigma_{n-1}=\sin^{n-2}\phi\,d\phi\,p^*d\sigma_{n-2}+\text{terms involving }\omega_{1s}\text{ or }\Omega^M_{\alpha\beta}.
\een
(See \er{wo Ome} for the definition of the curvature forms $\Omega^M_{\a\beta}$. Also compare \er{total Phi_k} when $k=0$ 
in view of \er{factorial}. In the case of one fixed sphere and its equator, \er{vol relation} without the extra terms is easy and follows from using spherical coordinates.) 
In the limit when $r\to 0$, the integrals of the terms involving $\omega_{1s}$ or $\Omega^M_{\a\beta}$ are zero by the same reason as 
in the last step of \er{start}. 

Therefore, continuing \er{start} and using iterated integrals, we have  
\begin{align*}
&\int_{\vn(M)}\Phi-\int_{\alpha_V(M)}\Phi=-\lim_{r\to 0}\int_{W_r} \wt{\Phi_0}=-\frac{1}{c_{n-1}}\lim_{r\to 0}\int_{W_r} d\sigma_{n-1}\\
=&-\frac{1}{c_{n-1}}\lim_{r\to 0}\int_{\alpha_{\partial V}(S_r^M(\sing\pmv))} \left(\int_{\phi(\alpha_V(x))}^0 \sin^{n-2}\phi\,d\phi\right) d\sigma_{n-2}\\
\overset{(1)}{=}&\frac{1}{c_{n-1}}\left(\int_{0}^\pi \sin^{n-2}\phi\,d\phi\right)
\lim_{r\to 0}\int_{\alpha_{\partial V}(S_r^M(\sing\pmv))} d\sigma_{n-2}\\
\overset{(2)}{=}&\frac{1}{c_{n-2}}\lim_{r\to 0}\int_{\alpha_{\partial V}(S_r^M(\sing\pmv))} d\sigma_{n-2}\overset{(3)}{=}\ind \pmv.
\end{align*}
Here equality (1) uses
\ben
\label{lim of phi}
\phi(\alpha_V(x))\to \pi \text{ for }x\in S_r^M(\sing\pmv),\text{ as }r\to 0,
\een
equality (2) uses the basic knowledge 
\ben\label{basic-know}
c_{n-1}= {c_{n-2}}\int_{0}^\pi \sin^{n-2}\phi\,d\phi,
\een 
and equality (3) is by the definition of index.
\end{proof}

\begin{remark} The construction of the chain $\cup_{x\in U_r}C_x$ is reminiscent of the topological method \cite{morse} of attaching $M\times I$ to $X$ and extending $V|_M$ to a vector field on $M\times I$ whose value at $(x,t)\in M\times I$ is $(1-t)V(x)+t\vec n(x)$. 
\end{remark}

\begin{remark} The homology group $H_{n-1}(STX|_M,\bz)\cong \bz\oplus \bz$ has two generators as  the image $\vec n(M)$ and a fiber sphere $ST_xM$ for $x\in M$ (see \cite{nie2}). Our proof shows that as a homology class, 
$$
\a_V(M)=\vec n(M)+(\ind \pmv) ST_xM.
$$
\end{remark}

\begin{remark} This first proof can also be presented very efficiently using the Mathai-Quillen form $\Psi$ \cite[III d)]{BZ} instead of our secondary Chern-Euler form $\Phi$. The author thanks a previous anonymous referee for pointing this out. 
\end{remark}

\section{Second proof by transgressing $\Phi$}
\label{trans}

In this section, we present a  study of the secondary Chern-Euler form $\Phi$ on $CSTM\subset STX|_M$ \er{eq-def-cstm}, leading to a proof of Theorem \ref{new main} and a second proof of Theorem \ref{two terms} using that.  

Recall the definition of the angle coordinate $\phi$ \er{def-phi}. Then $d\phi$ and $\dph$ are well-defined 1-form and vector field on $CSTM$. We write $d$ for exterior differentiation on $CSTM$, and $\iota_{\dph}$ for interior product with $\dph$. 

\begin{proposition}\label{much quicker} On $CSTM$, let 
\ben\label{write out}
\up=\iota_{\frac{\partial}{\partial \phi}}\Phi.
\een
Then the Lie derivative  
\ben
\label{new one}
\calL_\dph \Phi=d\up.
\een
Therefore
\ben\label{formal one}
\Phi-\pi^*\vec n^*\Phi=d\int_0^\phi \up \, dt.
\een
\end{proposition}

\begin{proof} \er{new one} follows from the Cartan homotopy formula (see, e.g., \cite[Prop. I.3.10]{KN})
$$
\calL_\dph \Phi=(d\,\iota_{\dph}+\iota_{\dph}d)\Phi=d\up,
$$
by \er{write out} and $d\Phi=0$ \er{d=0}. 

\er{formal one} then follows by integration since $\pi^*\vec n^*\Phi$ corresponds to the evaluation of $\Phi$ at $\phi=0$ by the definition of $\phi$ \er{def-phi} and we have for any fixed $\phi$
$$
\Phi-\pi^*\vec n^*\Phi=\int_0^\phi \calL_\dph \Phi\, dt=\int_0^\phi d\up \, dt=d\int_0^\phi \up\, dt.
$$
\end{proof}

Now we calculate $\up$ explicitly. Since $\Phi$ \er{phi} is invariant under $\so(n)$-changes of local frames, we adapt an idea from \cite{chern2} to use a nice frame for $TX|_M$ to facilitate the calculations about $\Phi$ on $CSTM$. 
Choose $e_1$ as in \ref{setup}. For $v\in CSTM$, let 
\ben
\label{en}
e_n=p(v)
\een
as defined in \er{def p}. Choose $e_2,\cdots,e_{n-1}$ so that $\{e_1,e_2,\cdots,e_{n-1},e_n\}$ is a positively oriented frame for $TX|_M$. (Therefore we need $n\geq 3$ from now on, with the $n=2$ case being simple.) Then in view of \er{def-phi}
\begin{equation}
\label{trig}
v=\cos\phi\, e_1+\sin\phi\, e_n.
\end{equation}


Let $(\Ome_{st})$ denote the curvature forms on $M$ of the induced metric from $X$. In view of \er{def-O}, 
\begin{gather}
\label{wo Ome}
\Ome_{st}=d\om_{st}-\sum_{r=2}^n \om_{sr}\om_{rt},
\end{gather}
\begin{gather}
\label{Ome}
\Om_{st}=\Ome_{st}+\om_{1s}\om_{1t}.
\end{gather}
Define the following differential forms on $STM$, regarded to be pulled back to $CSTM$ by $p$ \er{def p}, of degree $n-2$:
\begin{align}
\Phi^M(i,j)=\sum_\alpha \epsilon(\a) & \omega_{1\a_2}\cdots\omega_{1\a_{n-2i-j-1}}\Omega^M_{\a_{n-2i-j}\a_{n-2i-j+1}}\cdots \Omega^M_{\a_{n-j-2}\a_{n-j-1}}\label{phi ji}\\
&\omega_{\a_{n-j}n}\cdots\omega_{\a_{n-1}n}\nm
\end{align}
where the summations run over all permutations $\a$ of $\{2,\cdots,n-1\}$. It is easy to check that these $\Phi^M(i,j)$ are invariant under $\so(n-2)$-changes of the partial frames $\{e_2,\cdots,e_{n-1}\}$. Here the two parameters $i$ and $j$ stand for the numbers of curvature forms and $\om_{\alpha n}$'s involved. Define the following region of the indices $i,j$ 
\begin{align*}
D_1=\{(i,j)\in \bz\times\bz\, |\, i\geq 0,j\geq 0, 2i+j\leq n-2\}\label{d1}
\end{align*}
Then
$$
\Phi^M(i,j)\neq 0\Rightarrow (i,j)\in D_1. 
$$

\begin{remark} Our choice of the letter $\Phi^M$ 
is due to the following special case when there are no $\om_{1\a}$'s:
\begin{gather*}
\Phi^M(i,n-2i-2)=\Phi_i^M,
\end{gather*}
where $\Phi_i^M$ 
are forms on $STM$ defined by Chern \cite{chern2}. Since we are considering the case of boundary, we have the extra $\om_{1\a}$'s in our more general forms. 

Also note that the $\om_{1s}=0$ if the metric on $X$ is locally product near the boundary $M$. Therefore a lot of our forms vanish in that simpler case as considered in \cite{nie}. 
\end{remark}

We also introduce the following functions of $\phi$ \er{def-phi}, for non-negative integers $p$ and $q$, 
\begin{gather}
T(p,q)(\phi)=\cos^p\phi\sin^q\phi\label{def-t}\nm,\\
I(p,q)(\phi)=\int_0^\phi T(p,q)(t)\, dt\label{def-i}.
\end{gather}






\begin{proposition}\label{big Phi} We have the following concrete formulas
\begin{align}
\up=&\iota_\dph\Phi=\frac{1}{(n-2)!!c_{n-1}}\sum_{(i,j)\in D_1} a(i,j)(\phi)\, \Phi^M(i,j)\label{form of up}\\
a(i,j)(\phi)=&\sum_{k=i}^{[\frac{n-j}{2}]-1} (-1)^{n+j+k}\frac{(n-2k-2)!!}{2^k j!(n-2k-j-2)!i!(k-i)!}T(n-2k-j-2,j)(\phi)\label{aij}
\end{align}
\end{proposition}

\begin{proof} From \er{def-u}, \er{def-theta} and \er{trig}, we have
\begin{gather}
u_1=\cos\phi,\ u_n=\sin\phi,\ u_\alpha=0;\label{the u's}\\
\theta_1=-\sin\phi\, (d\phi+\om_{1n}),\ \theta_n=\cos\phi\,(d\phi+\om_{1n}),\label{theta 1 n}\\
\theta_\a=\cos\phi\, \om_{1\a}-\sin\phi\, \om_{\a n}\label{theta alpha}.
\end{gather}

From \er{the u's}, there are only two non-zero coordinates $u_1$ and $u_n$. Hence there are four cases for the positions of the indices $1$  and $n$ in $\Phi_k$ \er{eq-phi-j}:
\begin{enumerate}
\item $n-2k-1$ possibilities of $u_1\theta_n$
\item $2k$ possibilities of $u_1\Omega_{\a n}$
\item $n-2k-1$ possibilities of $u_n\theta_1$
\item $2k$ possibilities of $u_n\Omega_{1\a}$
\end{enumerate}
Only cases (i) and (iii) contribute $d\phi$ in view of \er{theta 1 n}, and hence we are only concerned with these two cases for the computation of $\up=\iota_\dph\Phi$. 
Starting with \er{eq-phi-j}, taking signs into considerations, by \er{the u's} and \er{theta 1 n}, by $\cos^2\phi+\sin^2\phi=1$, \er{theta alpha}, \er{Ome} and the multinomial theorem, we have
\begin{align}
\Phi_k=& (n-2k-1) (-1)^n \cos^2\phi  (d\phi+\omega_{1n})\label{total Phi_k}\\
&\quad \sum_\a \ep(\a)\theta_{\a_2}\cdots\theta_{\a_{n-2k-1}}\Omega_{\a_{n-2k}\a_{n-2k+1}}\cdots\Omega_{\a_{n-2}\a_{n-1}}\nm\\
&+(n-2k-1)(-1)^n \sin^2\phi  (d\phi+\omega_{1n})\nm\\
&\quad \sum_\a \ep(\a)\theta_{\a_2}\cdots\theta_{\a_{n-2k-1}}\Omega_{\a_{n-2k}\a_{n-2k+1}}\cdots\Omega_{\a_{n-2}\a_{n-1}}\nm\\
&+\cdots\nm\\
=&(-1)^n  (n-2k-1) (d\phi+\omega_{1n}) \sum_\a \ep(\a)\nm\\
&\quad (\cos\phi\,\om_{1\a_2}-\sin\phi\,\om_{a_2 n})\cdots (\cos\phi\,\om_{1\a_{n-2k-1}}-\sin\phi\,\om_{a_{n-2k-1} n})\nm\\
&\quad (\Omega^M_{\a_{n-2k}\a_{n-2k+1}}+\om_{1\a_{n-2k}}\om_{1\a_{n-2k+1}})\cdots(\Omega^M_{\a_{n-2}\a_{n-1}}+\om_{1\a_{n-2}}\om_{1\a_{n-1}})\nm\\
&+\cdots\nm\\
=&(-1)^n  (n-2k-1) (d\phi+\omega_{1n}) \nm\\
&\quad \mathop{\sum_{0\leq i\leq k}}_{0\leq j\leq n-2k-2}  \frac{(n-2k-2)!}{j!(n-2k-j-2)!} \cos^{n-2k-j-2}\phi\,(-\sin \phi)^j \frac{k!}{i!(k-i)!} \Phi^M(i,j)\nm\\
&+\cdots\nm\\
=&\mathop{\sum_{0\leq i\leq k}}_{0\leq j\leq n-2k-2} (-1)^{n+j}\frac{(n-2k-1)!k!}{j!(n-2k-j-2)!i!(k-i)!} T(n-2k-j-2,j)(\phi)\nm\\
&\qquad\qquad\quad (d\phi+\om_{1n})\Phi^M(i,j)+\cdots.\nm
\end{align}
From \er{phi} and the above, we get \er{form of up} and the coefficients $a(i,j)(\phi)$ in \er{aij}, 
after some immediate cancellations.
\end{proof}

\begin{definition}  For $(i,j)\in D_1$, define the following functions on $CSTM$
\begin{align}
&A(i,j)(\phi)=\int_0^\phi a(i,j)(t)\, dt\nm\\
=&\sum_{k=i}^{[\frac{n-j}{2}]-1} (-1)^{n+j+k}\frac{(n-2k-2)!!}{2^k j!(n-2k-j-2)!i!(k-i)!}I(n-2k-j-2,j)(\phi),\label{wo A}
\end{align}
in view of \er{aij} and \er{def-i}. Also define the differential form of degree $n-2$ on $CSTM$
\ben
\label{def Gamma}
\Gamma=\frac{1}{(n-2)!!c_{n-1}}\sum_{(i,j)\in D_1}A(i,j)(\phi)\Phi^M(i,j).
\een
\end{definition}

\begin{proof}[Proof of Theorem \ref{new main}] We just need to notice that $\Gamma=\int_0^\phi \up\, dt$ by Proposition \ref{big Phi} and use \er{formal one} in Proposition \ref{much quicker}.
\end{proof}


\begin{remark} Our first proof of Theorem \ref{new main} was through very explicit differentiations. Write 
$\Phi=d\phi\, \Upsilon+\Xi$ in view of \er{write out}. We can compute $\Xi$ explicitly. After correctly guessing the $\Gamma$ in \er{def Gamma}, we prove Theorem \ref{new main} by some differentiation formulas of differential forms in the spirit of \cite{chern2}, and some induction formulas for the functions $I(p,q)(\phi)$ in \er{def-i} through integration by parts.
\end{remark}

We finally arrive at
\begin{proof}[Second proof of Theorem \ref{two terms}] Let $B_r^M(\sing \partial V)$ (resp. $S_r^M(\sing \partial V))$ denote the union of smalls open balls (resp. spheres) of radii $r$ in $M$ around the finite set of points $\sing \partial V$. 
Then by $\partial V(x)=0\Leftrightarrow \a_V(x)=\pm \vec n(x)$, 
$$\a_V(M\backslash B^M_r(\sing \partial V))\subset CSTM.$$ 
By Theorem \ref{new main} and Stokes' theorem, 
\begin{align}
\label{stokes 1}
&\int_{\a_V(M)}\Phi-\int_{\vec n(M)}\Phi=\int_{\alpha_V(M)}\Phi-\pi^*\vec n^*\Phi=\lim_{r\to 0} \int_{\a_V(M\backslash B_r^M(\sing\partial V))}\Phi-\pi^*\vec n^*\Phi\\
=& \lim_{r\to 0} \int_{\a_V(M\backslash B_r^M(\sing\partial V))} d\Gamma
=-\lim_{r\to 0} \int_{\a_V(S_r^M(\sing \partial V))} \Gamma\nm\\
=&-\lim_{r\to 0} \int_{\a_V(S_r^M(\sing\partial V))}  \frac{1}{(n-2)!!c_{n-1}} A(0,n-2)(\phi)\Phi^M(0,n-2),\nm
\end{align}
since all the other $A(i,j)(\phi)\Phi^M(i,j)$ in \er{def Gamma}, for $(i,j)\in D_1$ and not equal to $(0,n-2)$, 
involve either curvature forms $\Omega^M_{\a\beta}$ or connection forms $\om_{1\a}$ and hence don't contribute in the limit when integrated over small spheres.

We have by \er{wo A} and \er{phi ji}
\begin{align}
&\frac{1}{(n-2)!!c_{n-1}} A(0,n-2)(\phi)\Phi^M(0,n-2)\label{vol again}\\
=&\frac{1}{(n-2)!!c_{n-1}}\frac{(n-2)!!}{(n-2)!} I(0,n-2)(\phi)\sum_\a \epsilon(\a) \om_{\a_2 n}\cdots\om_{\a_{n-1} n}\nm\\
=&\frac 1 {c_{n-1}} I(0,n-2)(\phi)\,p^*d\sigma_{n-2}\nm
\end{align}
with $\vol_{n-2}$ being the relative volume form of $S^{n-2}\to STM\to M$, since 
$$\sum_\a \epsilon(\a) \om_{\a_2 n}\cdots\om_{\a_{n-1} n}=(n-2)!\,p^*\vol_{n-2}$$
 in view of \er{en} and by comparison with \er{factorial}.

Continuing \er{stokes 1} and using \er{vol again}, we have
{\allowdisplaybreaks
\begin{align*}
& \int_{\alpha_V(M)}\Phi-\int_{\vec n(M)}\Phi\\
=&-\frac{1}{c_{n-1}} \lim_{r\to 0} \int_{\a_V(S_r^M(\sing\ppv)\cup S_r^M(\sing \pmv))} I(0,n-2)(\phi)\,p^*d\sigma_{n-2}\\
\overset{(1)}=& -\frac{1}{c_{n-1}} \bigl[I(0,n-2)(0)\lim_{r\to 0} \int_{\a_{\partial V}(S_r^M(\sing\ppv))} \vol_{n-2}\\
&\qquad\  \,+ I(0,n-2)(\pi)\lim_{r\to 0} \int_{\a_{\partial V}(S_r^M(\sing\pmv))} \vol_{n-2}\bigr]\\
\overset{(2)}=&\frac{1}{c_{n-2}} \lim_{r\to 0} \int_{\a_{\partial V}(S_r^M(\sing\pmv))} \vol_{n-2}\\
\overset{(3)}=&-\ind\pmv 
\end{align*}
}
Here equality (1) uses \er{pa}, \er{lim of phi}  and the similar
\begin{gather*}
\phi(\alpha_V(x))\to 0 \text{ for }x\in S_r^M(\sing\ppv),\text{ as }r\to 0.
\end{gather*}
In view of \er{def-i}, we have
\begin{gather*}
I(0,n-2)(0)=0,\ I(0,n-2)(\pi)=\int_0^\pi \sin^{n-2}\phi\,d\phi.
\end{gather*}
Then equality (2) follows from \er{basic-know}.
Equality (3) is by the definition of index. 
\end{proof}

\bigskip

\end{document}